\documentclass[12pt]{amsart}
\newcommand{\conservepaper}{
 \hoffset=-0.75in
 \setlength{\textwidth}{6.5in}
 \voffset=-0.5in
 \setlength{\textheight}{9.0in} 
 }
 \conservepaper

\usepackage{amssymb}

\newcommand{\Lie}[1]{\mathfrak{\lowercase{#1}}}
\newcommand{\real}{\mathord{\mathbb{R}}}
\newcommand{\AdG}{\operatorname{Ad}_G}
\newcommand{\adg}{\operatorname{ad}_{\Lie G}}
\newcommand{\zar}[1]{\overline{#1}}

\newcommand{\Rrank}{\operatorname{\hbox{$\real$-rank}}}
\newcommand{\unip}{\operatorname{unip}}
\newcommand{\iso}{\cong}
\newcommand{\Sl}{\operatorname{\Lie{SL}}}
\newcommand{\GL}{\operatorname{GL}}
\newcommand{\complex}{\mathord{\mathbb{C}}}
\newcommand{\SU}{\operatorname{SU}}
\newcommand{\SO}{\operatorname{SO}}
\newcommand{\so}{\operatorname{\Lie{so}}}
\newcommand{\Ortho}{\operatorname{O}}
\newcommand{\proj}{\operatorname{proj}}

 \swapnumbers

\numberwithin{equation}{section}
\newtheorem{prop}[equation]{Proposition}
\newtheorem{thm}[equation]{Theorem}
\newtheorem{lem}[equation]{Lemma}
\newtheorem{cor}[equation]{Corollary}

\theoremstyle{definition}
\newtheorem{defn}[equation]{Definition}
\newtheorem{notation}[equation]{Notation}
\newtheorem{ack}[equation]{Acknowledgments}

\theoremstyle{remark}
\newtheorem{rem}[equation]{Remark}

 \newcounter{case}
 \newenvironment{case}[1][\unskip]{\refstepcounter{case}\em
 \medskip \noindent {\bf Case \thecase\ #1.}\ }{\unskip\upshape}
 \renewcommand{\thecase}{\arabic{case}}

 \newcounter{subcase}
 \newenvironment{subcase}[1][\unskip]{\refstepcounter{subcase}\em
 \medskip \noindent {\bf Subcase \thesubcase\ #1.}\ }{\unskip\upshape}
\numberwithin{subcase}{case}

 \newcounter{subsubcase}
 \newenvironment{subsubcase}[1][\unskip]{\refstepcounter{subsubcase}\em
 \medskip \noindent {\bf Subsubcase \thesubsubcase\ #1.}\ }{\unskip\upshape}
\numberwithin{subsubcase}{subcase}

 \newcounter{subsubsubcase}
 \newenvironment{subsubsubcase}[1][\unskip]{\refstepcounter{subsubsubcase}\em
 \medskip \noindent {\bf Subsubsubcase \thesubsubsubcase\ #1.}\
}{\unskip\upshape} \numberwithin{subsubsubcase}{subsubcase}

\newenvironment{thmref}{\thmrefer}{}
\newcommand{\thmrefer}[1]{\renewcommand\theequation
 {\protect\ref{#1}$'$}\addtocounter{equation}{-1}}

\newcommand{\pref}[1]{{\upshape(\ref{#1})}}
\newcommand{\fullref}[2]{{\upshape\ref{#1}{\pref{#1-#2}}}}

\begin{document}

\title[Homogeneous Lorentz manifolds]
 {Homogeneous Lorentz manifolds \\ with simple isometry group}

\author{Dave Witte}

\address{Department of Mathematics, Oklahoma State University,
Stillwater, OK 74078}

\email{dwitte@math.okstate.edu,
http://www.math.okstate.edu/$\sim$dwitte}

\thanks{Submitted to \emph{Beitr\"age zur Algebra und Geometrie} (July 2000).
This version: 24~July 2000.} 

\begin{abstract}
 Let $H$ be a closed, noncompact subgroup of a simple Lie group~$G$, such that
$G/H$ admits an invariant Lorentz metric. We show that if $G = \SO(2,n)$, with
$n \ge 3$, then the identity component~$H^\circ$ of~$H$ is conjugate to
$\SO(1,n)^\circ$. Also, if $G = \SO(1,n)$, with $n \ge 3$, then $H^\circ$ is
conjugate to $\SO(1,n-1)^\circ$.
 \end{abstract}

\maketitle

\section{Introduction} \label{sect-intro}

\begin{defn}
 \begin{itemize}
 \item A \emph{Minkowski form} on a real vector space~$V$ is a nondegenerate
quadratic form that is isometric to the form $-x_1^2 + x_2^2 + \cdots +
x_{n+1}^2$ on~$\real^{n+1}$, where $\dim V = n+1 \ge 2$. 
 \item A \emph{Lorentz metric} on a smooth manifold~$M$ is a
choice of Minkowski metric on the tangent space~$T_p M$, for each $p \in M$,
such that the form varies smoothly as $p$ varies.
 \end{itemize}
 \end{defn}

A.~Zeghib \cite{ZeghibCpct} classified the compact homogeneous spaces that
admit an invariant Lorentz metric. In this note, we remove the assumption
of compactness, but add the restriction that the transitive group~$G$ is
almost simple. Our starting point is a special case of a theorem of N.~Kowalsky.

\begin{thm}[{N.~Kowalsky, cf.\ \cite[Thm.~5.1]{KowalskyAnnals}}]
\label{KowalskyTrans}
 Let $G/H$ be a nontrivial homogeneous space of a connected, almost simple Lie
group~$G$ with finite center. If there is a $G$-invariant Lorentz metric
on~$G/H$, then either
 \begin{enumerate}
 \item \label{KowalskyTrans-Riem}
 there is also a $G$-invariant Riemannian metric on $G/H$; or
 \item \label{KowalskyTrans-SO}
 $G$ is locally isomorphic to either $\SO(1,n)$ or $\SO(2,n)$, for
some~$n$.
 \end{enumerate}
 \end{thm}

As explained in the following elementary proposition, it is easy to characterize
the homogeneous spaces that arise in Conclusion~\pref{KowalskyTrans-Riem} of
Theorem~\ref{KowalskyTrans}, although it is probably not reasonable to expect a
complete classification.

\begin{notation}
 We use $\Lie G$ to denote the Lie algebra of a Lie group~$G$, and $\Lie H
\subset \Lie G$ to denote the Lie algebra of a Lie subgroup~$H$ of~$G$.
 \end{notation}

\begin{prop}[{cf.\ \cite[Thm.~1.1]{KowalskyAnnals}}]
 Let $G/H$ be a homogeneous space of a Lie group~$G$, such that $\Lie G$ is
simple and $\dim G/H \ge 2$. The following are equivalent.
 \begin{enumerate}
 \item The homogeneous space $G/H$ admits both a $G$-invariant Riemannian metric
and a $G$-invariant Lorentz metric.
 \item The closure of $\AdG H$ is compact, and leaves invariant a
one-dimensional subspace of~$\Lie G$ that is not contained in~$\Lie H$.
 \end{enumerate}
 \end{prop}

 The two main results of this note examine the cases that arise in
Conclusion~\pref{KowalskyTrans-SO} of Theorem~\ref{KowalskyTrans}.  It is well
known \cite[Egs.~2 and~3]{KowalskyCR} that $\SO(1,n)^\circ/\SO(1,n-1)^\circ$ and
$\SO(2,n)^\circ/\SO(1,n)^\circ$ have invariant Lorentz metrics. Also, for any
discrete subgroup~$\Gamma$ of $\SO(1,2)$, the Killing form provides an invariant
Lorentz metric on $\SO(1,2)^\circ/\Gamma$. We show that these are essentially
the only examples.

 Note that $\SO(1,1)$ and $\SO(2,2)$ fail to be almost simple. Thus, in
\fullref{KowalskyTrans}{SO}, we may assume
 \begin{itemize}
 \item $G$ is locally isomorphic to $\SO(1,n)$, and $n \ge 2$; or
 \item $G$ is locally isomorphic to $\SO(2,n)$, and $n \ge 3$.
 \end{itemize}

\begin{thmref}{so1n/H}
 \begin{prop} \label{SO1n/H}
  Let $G$ be a Lie group that is locally isomorphic to $\SO(1,n)$, with $n
\ge 2$. If $H$ is a closed subgroup of~$G$, such that 
 \begin{itemize}
 \item the closure of $\AdG H$ is not compact, and
 \item there is a $G$-invariant Lorentz metric on $G/H$,
 \end{itemize}
 then either
 \begin{enumerate}
 \item after any identification of~$\Lie G$ with $\so(1,n)$, the
subalgebra $\Lie H$ is conjugate to a standard copy of $\so(1,n-1)$ in
$\so(1,n)$, or
 \item $n = 2$ and $H$ is discrete.
 \end{enumerate}
  \end{prop}
 \end{thmref}

\begin{thmref}{so2n/H}
 \begin{thm}  \label{SO2n/H}
  Let $G$ be a Lie group that is locally isomorphic to $\SO(2,n)$, with $n
\ge 3$. If $H$ is a closed subgroup of~$G$, such that 
 \begin{itemize}
 \item the closure of $\AdG H$ is not compact, and
 \item there is a $G$-invariant Lorentz metric on $G/H$,
 \end{itemize}
 then, after any identification of~$\Lie G$ with $\so(2,n)$, the subalgebra
$\Lie H$ is conjugate to a standard copy of $\so(1,n)$ in
$\so(2,n)$.
  \end{thm}
 \end{thmref}

  N.~Kowalsky announced a much more general result than Theorem~\ref{SO2n/H} in
\cite[Thm.~4]{KowalskyCR}, but it seems that she did not publish a proof before
her premature death. She announced a version of Proposition~\ref{SO1n/H} (with
much more general hypotheses and a somewhat weaker conclusion) in
\cite[Thm.~3]{KowalskyCR}, and a proof appears in her Ph.D.\ thesis
\cite[Cor.~6.2]{KowalskyThesis}.

\begin{rem}
 It is easy to see that there is a $G$-invariant Lorentz metric on $G/H$ if and
only if there is an $(\AdG H)$-invariant Minkowski form on $\Lie G/ \Lie H$.
Thus, although Proposition~\ref{SO1n/H} and Theorem~\ref{SO2n/H} are
geometric in nature, they can be restated in more algebraic terms. It is in
such a form that they are proved in \S\ref{so1n-sect} and~\S\ref{so2n-sect}.
 \end{rem}

Proposition~\ref{SO1n/H} and Theorem~\ref{SO2n/H} are used in work of S.~Adams
\cite{AdamsNontame} on nontame actions on Lorentz manifolds. See
\cite{ZimmerLorentz, KowalskyAnnals, AdamsStuckLorentz, ZeghibId, AdamsLocfree,
AdamsTrans} for some other research concerning actions of Lie groups on Lorentz
manifolds.

\begin{ack}
 The author would like to thank the Isaac Newton Institute for Mathematical
Sciences for providing the stimulating environment where this work was
carried out. It is also a pleasure to thank Scot Adams for suggesting this
problem and providing historical background.
 The research was partially supported by a grant from the National Science
Foundation (DMS-9801136).
 \end{ack}

\section{Homogeneous spaces of $\SO(1,n)$} \label{so1n-sect}

The following lemma is elementary.

\begin{lem} \label{lem-std-rep}
  Let $\pi$ be the standard representation of $\Lie G = \so(1,k)$ on
$\real^{k+1}$, and let $\Lie G = \Lie K + \Lie A + \Lie N$ be an Iwasawa
decomposition of~$\Lie G$.
  \begin{enumerate}
  \item
  The representation~$\pi$ has only one positive weight {\upshape(}with respect
to~$\Lie A${\upshape)}, and the corresponding weight space is 1-dimensional.
  \item
  There are subspaces $V$ and~$W$ of~$\real^{k+1}$, such that
  \begin{enumerate}
  \item
  $\dim ( \real^{k+1}/V) = 1$;
  \item
  $\dim W = 1$;
  \item
  $\pi( \Lie N ) V \subset W$;
  \item
  for all nonzero $u \in \Lie N$, we have $\pi(u)^2 \real^{k+1} = W$; and
  \item
  for all nonzero $u \in \Lie N$ and $v \in \real^{k+1}$, we have $\pi(u)^2 v =
0$ if and only if $v \in V$.
  \end{enumerate}
  \end{enumerate}
  \end{lem}

\begin{cor}  \label{so1k}
  Let $\Lie H$ be a subalgebra of a real Lie algebra~$\Lie G$, let $Q$ be
a Minkowski form on $\Lie G/\Lie H$, and define
  $\pi \colon N_G(\Lie H) \to \GL(\Lie G/\Lie H)$ by
  $\pi(g)(v + \Lie H) = (\AdG g)v + \Lie H$.
  \begin{enumerate}
  \item \label{so1k-A}
  Suppose $T$ is a connected Lie subgroup of~$G$ that normalizes~$H$, such
that $\pi(T) \subset \SO(Q)$ and $\AdG T$ is diagonalizable
over~$\real$. Then, for any ordering of the $T$-weights on~$\Lie G$, the
subalgebra $\Lie H$ contains codimension-one subspaces of both $\Lie G^+$ and
$\Lie G^-$, where $\Lie G^+$ is the sum of all the
positive weight spaces of~$T$, and $\Lie G^-$ is the sum of all the negative
weight spaces of~$T$.
  
  \item \label{so1k-N}
   If $U$ is a connected Lie subgroup of~$G$ that normalizes~$H$, such that
$\pi(U) \subset \SO(Q)$ and $\AdG U$ is unipotent, then there are
subspaces $V/ \Lie H$ and~$W / \Lie H$ of~$\Lie G / \Lie H$, such that
  \begin{enumerate}
  \item \label{so1k-N-dimV}
  $\dim ( \Lie G / V) = 1$;
  \item \label{so1k-N-dimW}
  $\dim (W / \Lie H) = 1$;
  \item \label{so1k-N-adVinW}
  $[V, \Lie U] \subset W$;
  \item \label{so1k-N-ad2inW}
  for each $u \in \Lie U$, either $W = \Lie H + (\adg u)^2 \Lie G$, or $[\Lie G,
u] \subset \Lie H$; and
  \item \label{so1k-N-ad2v=0}
  for all $u \in \Lie U$, we have $(\adg u)^2 V \subset \Lie H$.
  \end{enumerate}
  \end{enumerate}
  \end{cor}

\begin{prop} \label{so1n/H}
  Let $H$ be a Lie subgroup of $G = \SO(1,n)$, with $n \ge 2$, such that 
 \begin{itemize}
 \item the closure of~$H$ is not compact; and
 \item there is an $(\AdG H)$-invariant Minkowski form on $\Lie G / \Lie H$.
 \end{itemize} Then either
 \begin{enumerate}
 \item $H^\circ$ is conjugate to a standard copy of $\SO(1,n-1)^\circ$ in
$\SO(1,n)$, or
 \item \label{so1n/H-n=2}
 $n = 2$ and $H^\circ$ is trivial.
 \end{enumerate}
  \end{prop}

\begin{proof}
  Let $\zar H$ be the Zariski closure of~$H$, and note that the Minkowski form
is also invariant under~$\AdG {\zar H}$. Replacing $H$ by a finite-index
subgroup, we may assume $\zar H$ is Zariski connected.

Let $G = KAN$ be an Iwasawa decomposition of~$G$. 

\setcounter{case}{0}

\begin{case}
  Assume $n \ge 3$ and $A \subset \zar H$.
  \end{case}
  From Corollary~\fullref{so1k}{A}, we see that $\Lie H$ contains codimension-one
subspaces of both $\Lie N$ and~$\Lie N^-$. (Note that this implies $H^\circ$
is nontrivial.) This implies that $\zar H$ is reductive.
(Because $(H \cap N)^\circ \unip \zar H$ is a unipotent subgroup that
intersects~$N$ nontrivially (and $\Rrank G = 1$), it must be contained in~$N$,
so $\unip \zar H \subset N$. Similarly, $\unip \zar H \subset N^-$. Therefore
$\unip \zar H \subset N \cap N^- = e$.) Then, since $\zar H$ contains a
codimension-one subgroup of~$N$, and since $A\subset\zar H$, it follows that
$\zar H$ is conjugate to either $\SO(1,n-1)$ or $\SO(1,n)$. Because $H^\circ$ is
a nontrivial, connected, normal subgroup of~$\zar H$, we conclude that $H^\circ$
is conjugate to either $\SO(1,n-1)^\circ$ or $\SO(1,n)^\circ$. Because $\Lie G /
\Lie H \neq 0$ (else $\dim \Lie G / \Lie H = 0 < 2$, which contradicts the fact
that there is a Minkowski form on $\Lie G/\Lie H$), we see that $H^\circ$ is
conjugate to $\SO(1,n-1)^\circ$.

\begin{case}
  Assume $n \ge 3$ and $\zar H$ does not contain any nontrivial hyperbolic
elements.
  \end{case}
  The Levi subgroup of $\zar H$ must be compact, and the radical of $\zar H$
must be unipotent, so choose a compact $M$ and a nontrivial unipotent
subgroup $U$ such that $\zar H = M \ltimes U$.
Replacing $H$ by a conjugate,
we may assume, without loss of generality, that $U \subset N$.

  Let us show, for every nonzero $u \in \Lie U$, that $[\Lie G, u] \not\subset
\Lie H$. From the Morosov Lemma \cite[Thm.~17(1), p.~100]{Jacobson}, we know
there exists $v \in \Lie G$, such that $[v,u]$ is hyperbolic (and nonzero). If
$[v,u] \in \Lie H$, this contradicts the fact that $\zar H$ does not contain
nontrivial hyperbolic elements.

Let $V/\Lie H$ and $W/\Lie H$ be subspaces of $\Lie G / \Lie H$ as in
Corollary~\fullref{so1k}{N}. Because $(\adg u)^2 \Lie G = \Lie N$ for every
nonzero $u \in \Lie N$, we have $W = \Lie N + \Lie H$
(see~\fullref{so1k}{N-ad2inW}), so $\dim \Lie N / (\Lie H \cap \Lie N) = 1$
(see~\fullref{so1k}{N-dimW})
  and
  \begin{equation} \label{uVinN+M}
  [\Lie U, V] \subset W = \Lie N + \Lie H \subset \Lie N +
\Lie{\zar H} = \Lie N + \Lie M
  \end{equation}
  (see~\fullref{so1k}{N-adVinW}).

Assume, for the moment, that $n \ge 4$. Then
  \begin{eqnarray*}
  \dim \Lie U + \dim (V \cap \Lie N^-)
  &\ge& \dim ( \Lie H \cap \Lie N ) + \dim( V \cap \Lie N^- )
  \ge (\dim \Lie N - 1) + (\dim \Lie N^- - 1) \\
  &=& (n-2) + (n-2)
  \ge n
  > \dim \Lie N.
  \end{eqnarray*}
  This implies that there exist $u \in \Lie U$ and $v
\in V \cap \Lie N^-$, such that $\langle u,v \rangle \iso \Sl(2,\real)$, with
$[u,v]$ hyperbolic (and nonzero). This contradicts the fact that $\Lie M + \Lie
N$ has no nontrivial hyperbolic elements.

We may now assume that $n = 3$. For any nonzero $u \in \Lie N$, we have
  $$ \dim [u, V] \ge \dim [u, \Lie G] - 1
  = \dim \Lie N + 1 > \dim \Lie N ,$$
  so $[\Lie U, V] \not\subset \Lie N$. Then, from~(\ref{uVinN+M}), we conclude
that $\Lie M \neq 0$, so $\Lie M$ acts irreducibly on~$\Lie N$. This
contradicts the fact that $\Lie H \cap \Lie N$ is a codimension-one
subspace of~$\Lie N$ that is normalized by~$\Lie M$.

\begin{case}
 Assume $n = 2$.
 \end{case}
 We may assume $H^\circ$ is nontrivial (otherwise Conclusion~\pref{so1n/H-n=2}
holds). We must have $\dim \Lie G / \Lie H \ge 2$, so we conclude that $\dim
H^\circ = 1$ and $\dim \Lie G / \Lie H = 2$. Because $\SO(1,1)$ consists of
hyperbolic elements, this implies that $\AdG h$ acts diagonalizably on $\Lie G/
\Lie H$, for every $h \in H$. Therefore $H^\circ$ is conjugate to~$A$, and,
hence, to $\SO(1,1)^\circ$.
  \end{proof}

\section{Homogeneous spaces of $\SO(2,n)$} \label{so2n-sect}

\begin{thm}[{Borel-Tits \cite[Prop.~3.1]{BorelTitsunip}}]
 \label{Titsunip}
  Let $H$ be an $F$-subgroup of a reductive algebraic group~$G$ over a field~$F$
of characteristic zero. Then there is a parabolic $F$-subgroup $P$ of~$G$, such
that $\unip H \subset \unip P$ and $H \subset N_G(\unip H) \subset P$.
  \end{thm}

\begin{notation}
  Let $k = \lfloor n/2 \rfloor$.
Identifying $\complex^{k+1}$ with $\real^{2k+2}$ yields
an embedding of $\SU \bigl( 1, k \bigr)$ in $\SO(2,2k)$. Then the
inclusion $\real^{2k+2} \hookrightarrow \real^{2k+3}$ yields an
embedding of $\SU \bigl( 1, k \bigr)$ in $\SO(2, 2k+1)$. Thus, we may identify
$\SU \bigl( 1, \lfloor n/2 \rfloor \bigr)$ with a subgroup of $\SO(2,n)$.
  \end{notation}

  We use the following well-known result to shorten one case of the proof of
Theorem~\ref{so2n/H}.

\begin{lem}[{\cite[Lem.~6.8]{OhWitte-CK}}]
  If $L$ is a connected, almost-simple subgroup of $\SO(2,n)$, such that
$\Rrank L = 1$ and $\dim L > 3$, then $L$ is conjugate under
$\operatorname{O}(2,n)$ to a subgroup of either $\SO(1,n)$ or $\SU \bigl( 1,
\lfloor n/2 \rfloor \bigr)$.
  \end{lem}

\begin{cor} \label{dim(HcapA)}
  Let $L$ be a connected, reductive subgroup of $G = \SO(2,n)$, such 
that $\Rrank
L = 1$. Then $\dim U \le n-1$, for every connected, unipotent
subgroup~$U$ of~$L$.

Furthermore, if $\dim U = n-1$, then either
  \begin{enumerate}
  \item $L$ is conjugate to $\SO(1,n)^\circ$; or
  \item $n$ is even, and $L$ is conjugate under $\Ortho(2,n)$ to $\SU(1, n/2)$.
  \end{enumerate}
  \end{cor}

\begin{thm}  \label{so2n/H}
  Let $H$ be a Lie subgroup of $G = \SO(2,n)$, with $n \ge 3$, such that 
 \begin{itemize}
 \item the closure of~$H$ is not compact, and
 \item there is an $(\AdG H)$-invariant Minkowski form on $\Lie G / \Lie H$.
 \end{itemize}
 Then $H^\circ$ is conjugate to a standard copy of $\SO(1,n)^\circ$ in
$\SO(2,n)$.
  \end{thm}

\begin{proof}
  Let $\zar H$ be the Zariski closure of~$H$, and note that the Minkowski form
is also invariant under~$\AdG {\zar H}$. Replacing $H$ by a finite-index
subgroup, we may assume $\zar H$ is Zariski connected.

  Let $G = KAN$ be an Iwasawa decomposition of~$G$. For each real root~$\phi$
of~$\Lie G$ (with respect to the Cartan subalgebra~$\Lie A$), let $\Lie G_\phi$
be the corresponding root space, and let
  $\proj_\phi \colon \Lie G \to \Lie G_\phi$
  and
  $\proj_{\phi \oplus -\phi} \colon \Lie G \to \Lie G_\phi + \Lie G_{-\phi}$
  be the natural projections. Fix a choice of simple real roots $\alpha$
and~$\beta$ of~$\Lie G$, such that $\dim \Lie G_\alpha = 1$ and $\dim \Lie
G_\beta = n-2$ (so the positive real roots are $\alpha$, $\beta$,
$\alpha+\beta$, and $\alpha + 2\beta$). Replacing $N$ by a conjugate under the
Weyl group, we may assume $\Lie N = \Lie G_\alpha + \Lie G_\beta + \Lie
G_{\alpha+\beta} + \Lie G_{\alpha+2\beta}$. From the classification of
parabolic subgroups \cite[Prop.~5.14, p.~99]{BorelTitsreductive}, we know
that the only proper parabolic subalgebras of~$\Lie G$ that contain $\Lie
N_{\Lie G}(\Lie N)$ are
  \begin{equation} \label{parabs}
  \mbox{$\Lie N_{\Lie G}(\Lie N)$,
  \ $\Lie P_\alpha = \Lie N_{\Lie G}(\Lie N) + \Lie G_{-\alpha}$,
  \  and
  \ $\Lie P_\beta = \Lie N_{\Lie G}(\Lie N) + \Lie G_{-\beta}$.}
  \end{equation}

\setcounter{case}{0}

\begin{case}
  Assume $\Lie{\zar H}$ contains nontrivial hyperbolic elements.
  \end{case}
  Let $\Lie T = \Lie {\zar H} \cap \Lie A$. Replacing $H$ by a conjugate, we may
assume $\Lie T \neq 0$.

\begin{subcase}
  Assume $\Lie T \in \{ \ker (\alpha+\beta), \ker \beta \}$.
  \end{subcase}

\begin{subsubcase}
  Assume $\zar H$ is reductive.
  \end{subsubcase}
  We may assume $\Lie T = \ker (\alpha+\beta)$ (if necessary, replace $H$ with
its conjugate under the Weyl reflection corresponding to the root~$\alpha$).
Then, from Corollary~\fullref{so1k}{A}, we see that $\Lie H$ contains
a codimension-one subspace of $\Lie G_{\alpha+2\beta} + \Lie G_{\beta} + \Lie
G_{-\alpha}$. (Note that this implies $H^\circ$ is nontrivial.)

  Let $\Lie N' = \Lie G_{\alpha+\beta} + \Lie G_{\alpha+2\beta} + \Lie G_{\beta}
+ \Lie G_{-\alpha}$, so $\Lie N'$ is the Lie algebra of a maximal unipotent
subgroup of~$G$. (In fact, $\Lie N'$ is the image of~$\Lie N$ under the Weyl
reflection corresponding to the root~$\alpha$.) From the preceding 
paragraph, we
know that
  $$\dim (\Lie {\zar H} \cap \Lie N')
  \ge \dim (\Lie G_{\alpha+2\beta} + \Lie G_{\beta} + \Lie G_{-\alpha}) - 1
  = n-1
  .$$
  Therefore, Corollary~\ref{dim(HcapA)} implies that $\zar H$ is 
conjugate (under
$\Ortho(2,n)$) to either $\SO(1,n)$ or $\SU(1,n/2)$.  It is easy to see that
$\zar H$ is not conjugate to $\SU(1,n/2)$.  (See \cite[proof of
Thm.~1.5]{OhWitte-CK} for an explicit description of $\Lie{SU}(1,n/2) \cap
\Lie N$. If $n$ is even, then $n > 3$, so $\Lie{SU}(1,n/2)$ does not contain a
codimension-one subspace of any $(n-2)$-dimensional root space, but $\zar{\Lie
H}$ does contain a codimension-one subspace of $\Lie G_{\beta}$.) Therefore, we
conclude that $\zar H$ is conjugate to $\SO(1,n)$. Then, because
$H^\circ$ is a nontrivial, connected, normal subgroup of~$\zar H$, we
conclude that $H^\circ = (\zar H)^\circ$ is conjugate to
$\SO(1,n)^\circ$.

\begin{subsubcase}
   Assume $\zar H$ is not reductive.
  \end{subsubcase}
  Let $P$ be a maximal parabolic subgroup of~$G$ that contains~$\zar H$
(see Theorem~\ref{Titsunip}). By replacing $P$ and~$H$ with conjugate
subgroups, we may assume that $P$ contains the minimal parabolic
subgroup $N_G(N)$. Therefore, the classification of parabolic
subalgebras (\ref{parabs}) implies that $P$ is either $P_\alpha$
or~$P_\beta$.

\begin{subsubsubcase}
  Assume $\Lie T = \ker (\alpha+\beta)$.
  \end{subsubsubcase}
  From Corollary~\fullref{so1k}{A}, we see that $\Lie H$ (and hence also~$\Lie
P$) contains codimension-one subspaces of $\Lie G_{\alpha+2\beta} + \Lie
G_{\beta} + \Lie G_{-\alpha}$ and $\Lie G_{-\alpha-2\beta} + \Lie G_{-\beta} +
\Lie G_{\alpha}$. Because $\Lie P_\alpha$ does not
contain such a subspace of $\Lie G_{-\alpha-2\beta} + \Lie G_{-\beta} + \Lie
G_{\alpha}$, we conclude that $P = P_\beta$. Furthermore, because the
intersection of~$\Lie P_\beta$ with each of these subspaces does have
codimension~one, we conclude that $\Lie H$ has precisely the same intersection;
therefore
  $(\Lie G_{\alpha+2\beta} + \Lie G_{\beta}) + (\Lie G_{-\beta} + \Lie
G_{\alpha}) \subset \Lie H$. Hence
  $\Lie H \supset [\Lie G_\alpha, \Lie G_\beta] = \Lie G_{\alpha+\beta}$.
  We now have
  $$ (\adg  \Lie G_{\alpha+\beta})^2 \Lie G
  = \Lie G_\alpha + \Lie G_{\alpha+\beta} + \Lie G_{\alpha+2\beta}
  \equiv 0 \pmod{\Lie H} ,$$
  so Corollary~\fullref{so1k}{N-ad2inW} implies
  $$ \Lie H \supset [\Lie G, \Lie G_{\alpha+\beta}]
  \supset [\Lie G_{-\alpha-\beta}, \Lie G_{\alpha+\beta}]
  \supset \ker \beta .$$
  This contradicts the fact that $\Lie {\zar H} \cap \Lie A = \Lie T = \ker
(\alpha+\beta)$.

\begin{subsubsubcase}
  Assume $\Lie T = \ker \beta$.
  \end{subsubsubcase}
  From Corollary~\fullref{so1k}{A}, we see that $\Lie H$ (and hence also~$\Lie
P$) contains a codimension-one subspace of $\Lie G_{-\alpha} + \Lie
G_{-\alpha-\beta} + \Lie G_{-\alpha-2\beta}$. Because neither $\Lie P_\alpha$
nor  $\Lie P_\beta$ contains such a subspace, this is a contradiction.

\begin{subcase} \label{reductive-alpha}
  Assume $\Lie T \in \{ \ker \alpha, \ker(\alpha+ 2\beta) \}$.
  \end{subcase}
  We may assume $\Lie T = \ker \alpha$ (if necessary, replace $H$ with
its conjugate under the Weyl reflection corresponding to the root~$\beta$).
  From Corollary~\fullref{so1k}{A}, we see that $\Lie H$ contains a 
codimension-one
subspace of $\Lie G_{\beta} + \Lie G_{\alpha+\beta} + \Lie G_{\alpha+2\beta}$.
Because any codimension-one subalgebra of a nilpotent Lie algebra must contain
the commutator subalgebra, we conclude that $\Lie H$
contains $\Lie G_{\alpha+2\beta}$. Then we have
  $$ (\adg \Lie G_{\alpha+2\beta})^2 \Lie G = \Lie G_{\alpha+2\beta} \equiv 0
\pmod{\Lie H} ,$$
  so Corollary~\fullref{so1k}{N-ad2inW} implies
  $$\Lie H \supset [\Lie G, \Lie G_{\alpha+2\beta}]
  \supset \Lie G_{\beta} + \Lie G_{\alpha+\beta} + \Lie G_{\alpha+2\beta} .$$
  Similarly, we also have $\Lie H \supset \Lie G_{-\beta} + \Lie
G_{-\alpha-\beta} + \Lie G_{-\alpha-2\beta}$. It is now easy to show that $\Lie
H \supset \Lie G_\phi$ for every real root~$\phi$, so $\Lie H = \Lie G$. This
contradicts the fact that $\Lie G / \Lie H \neq 0$.

\begin{subcase}
  Assume $\Lie T$ contains a regular element of~$\Lie A$.
  \end{subcase}
  Replacing $H$ by a conjugate under the Weyl group, we may assume that $\Lie
N$ is the sum of the positive root spaces, with respect to~$\Lie T$. Then,
from Corollary~\fullref{so1k}{A}, we see that $\Lie H$ contains codimension-one
subspaces of both $\Lie N$ and~$\Lie N^-$. Therefore, $\Lie H$ contains
codimension-one subspaces of $\Lie G_{\beta} + \Lie G_{\alpha+\beta} + \Lie
G_{\alpha+2\beta}$ and $\Lie G_{-\beta} + \Lie G_{-\alpha-\beta} + \Lie
G_{-\alpha-2\beta}$, so the argument of Subcase~\ref{reductive-alpha} applies.

\begin{case} \label{nohyper}
  Assume $\Lie{\zar H}$ does not contain nontrivial hyperbolic elements.
  \end{case}
  The Levi subgroup of $\zar H$ must be compact, and the radical of $\zar H$
must be unipotent, so choose a compact $M$ and a nontrivial unipotent
subgroup $U$ such that $\zar H = M \ltimes U$.
Choose subspaces $V/\Lie H$ and $W/\Lie H$ of $\Lie G/ \Lie H$ as in
Corollary~\fullref{so1k}{N}.

Let $P$ be a proper parabolic subgroup of~$G$, such that $U \subset \unip P$
and $H \subset P$ (see Theorem~\ref{Titsunip}). Replacing $H$ and~$P$ by
conjugates, we may assume, without loss of generality, that $P$ contains the
minimal parabolic subgroup $N_G(N)$ (so $\unip P \subset N$).
  From the classification of parabolic subalgebras (\ref{parabs}), we know that
there are only three possibilities for~$P$. We consider each of these
possibilities separately.

  First, though, let us show that
  \begin{equation} \label{Ufaithful}
  \mbox{for every nonzero $u \in \Lie U$, we have
$[\Lie G, u] \not\subset \Lie H$.}
   \end{equation}
  From the Morosov Lemma \cite[Thm.~17(1), p.~100]{Jacobson}, we know there
exists $v \in \Lie G$, such that $[v,u]$ is hyperbolic (and nonzero). If $[v,u]
\in \Lie H$, this contradicts the fact that $\Lie{\zar H}$ does not contain
nontrivial hyperbolic elements.

\begin{subcase} \label{so2npf-Pminl}
  Assume $P = N_G(N)$ is a minimal parabolic subgroup of~$G$.
  \end{subcase}

\begin{subsubcase}
  Assume $\proj_{\beta} \Lie U \neq 0$.
  \end{subsubcase}
  Choose $u \in \Lie U$, such that $\proj_{\beta} u \neq 0$, and let $Z = (\adg
u)^2 \Lie G_{-\alpha-2\beta}$. (So $\dim Z = 1$, $\proj_{-\alpha} Z \neq 0$,
and $\proj_{-\alpha-\beta} Z = 0$.) From Corollary~\fullref{so1k}{N-ad2inW}, we
know that
  $Z \subset W$.
  Then, because
  $\proj_{-\alpha} \Lie H \subset \proj_{-\alpha} \Lie P = 0$,
  we conclude, from Corollary~\fullref{so1k}{N-dimW}, that
  $W = \Lie H + Z$.

Because $W = \Lie H + Z \subset \Lie P + Z$, we
have $\proj_{-\alpha-\beta} W = 0$. Therefore, because $\proj_\beta u \neq 0$,
we conclude, from Corollary~\fullref{so1k}{N-adVinW}, that
$\proj_{-\alpha-2\beta} V = 0$, so Corollary~\fullref{so1k}{N-dimV} implies
that $V = \ker (\proj_{-\alpha-2\beta})$. In particular, we have $\Lie
G_{-\beta} \subset V$, so Corollary~\fullref{so1k}{N-adVinW} implies $[\Lie
G_{-\beta}, u] \subset W$. Therefore, we have
  \begin{eqnarray*}
  [\Lie G_{-\beta}, \proj_\beta u]
  &\subset& [\Lie G_{-\beta}, u + (\Lie G_\alpha +
\Lie G_{\alpha+\beta} + \Lie
G_{\alpha+2\beta})]
  = [\Lie G_{-\beta}, u] + [\Lie G_{-\beta}, \Lie G_\alpha + \Lie
G_{\alpha+\beta} + \Lie G_{\alpha+2\beta}] \\
  &\subset& W + (\Lie G_\alpha + \Lie G_{\alpha+\beta})
  = \Lie H + Z + (\Lie G_\alpha + \Lie G_{\alpha+\beta})
  \subset \Lie M + \Lie N + Z .
  \end{eqnarray*}
  Because $\proj_{-\alpha} [\Lie G_{-\beta}, \proj_\beta u] = 0$, we conclude
that $[\Lie G_{-\beta}, \proj_\beta u] \subset \Lie M + \Lie N$.
  This contradicts the fact that $\Lie M + \Lie N$ does not contain nontrivial
hyperbolic elements.

\begin{subsubcase}
  Assume $\proj_{\beta} \Lie U = 0$.
  \end{subsubcase}
  Replacing $H$ by a conjugate under~$N$, we may assume $\Lie M \subset \Lie
G_0$, so $\proj_{\beta} \Lie {\zar H} = 0$.

  We have $\Lie U \subset \Lie G_\alpha + \Lie G_{\alpha+\beta} + \Lie
G_{\alpha+2\beta}$, so
  $(\adg u)^2 \Lie G \subset \Lie G_\alpha + \Lie G_{\alpha+\beta} + \Lie
G_{\alpha+2\beta}$ for every $u \in \Lie U$. Thus,
Corollary~\fullref{so1k}{N-ad2inW} implies
  $W \subset (\Lie G_\alpha + \Lie G_{\alpha+\beta} +
\Lie G_{\alpha+2\beta}) + \Lie H$.

We have
  $$\proj_{\beta \oplus -\beta} W
  \subset \proj_{\beta \oplus -\beta} (\Lie G_\alpha + \Lie G_{\alpha+\beta} +
\Lie G_{\alpha+2\beta}) + \proj_{\beta \oplus -\beta} \Lie H
  = 0 ,$$
  so Corollary~\fullref{so1k}{N-adVinW} implies that $\proj_{\beta \oplus
-\beta} \bigl( (\adg \Lie U) V \bigr) = 0$.

\begin{subsubsubcase} \label{projalphanot0}
  Assume $\proj_{\alpha} u \neq 0$, for some $u \in \Lie U$.
  \end{subsubsubcase}
  From the conclusion of the preceding paragraph, we know that $\proj_{-\beta}
\bigl( (\adg u) V \bigr) = 0$. Because $\proj_\beta u = 0$ and $\proj_\alpha
\neq 0$, this implies $\proj_{-\alpha-\beta} V = 0$, so $V = \ker
(\proj_{-\alpha-\beta})$ (see~\fullref{so1k}{N-dimV}). In particular, $\Lie
G_{-\alpha} \subset V$, so Corollary~\fullref{so1k}{N-adVinW} implies
  \begin{eqnarray*}
  [\Lie G_{\alpha}, \Lie G_{-\alpha}]
  &\subset& [\Lie U + (\Lie G_{\alpha+\beta} + \Lie
G_{\alpha+2\beta}) , \Lie G_{-\alpha}]
  \subset [\Lie U , V] + [\Lie G_{\alpha+\beta} + \Lie
G_{\alpha+2\beta}, \Lie G_{-\alpha}] \\
  &\subset& W + \Lie G_\beta
  \subset \Lie H + \Lie N
  \subset \Lie M + \Lie N .
  \end{eqnarray*}
  This contradicts the fact that $\Lie M + \Lie N$ does not contain nontrivial
hyperbolic elements.

\begin{subsubsubcase}
  Assume $\proj_{\alpha+\beta} u \neq 0$, for some $u \in \Lie U$.
  \end{subsubsubcase}
  From Subsubsubcase~\ref{projalphanot0}, we may assume $\proj_{\alpha} u = 0$.
  Because $0 = \proj_{\beta \oplus -\beta} \bigl( (\adg u) V \bigr)$ has
codimension $\le 1$ in
  $\proj_{\beta \oplus -\beta} \bigl( (\adg u) \Lie G \bigr)$
  (see~\fullref{so1k}{N-dimV}),
  which contains the 2-dimensional subspace
  $\proj_{\beta \oplus -\beta} \bigl( [u, \Lie G_{-\alpha-2\beta} + \Lie
G_{-\alpha}] \bigr)$,
  we have a contradiction.

\begin{subsubsubcase}
  Assume $\Lie U = \Lie G_{\alpha+2\beta}$.
  \end{subsubsubcase}
  (This argument is similar to Subsubsubcase \ref{projalphanot0}.)
  Because $\proj_\beta \bigl( (\adg \Lie U) V \bigr) = 0$, we know that
$\proj_{-\alpha-\beta} V = 0$, so $V = \ker (\proj_{-\alpha-\beta})$
(see \fullref{so1k}{N-dimV}). In particular, $\Lie G_{-\alpha-2\beta} \subset
V$, so Corollary~\fullref{so1k}{N-adVinW} implies
  $$[\Lie G_{\alpha+2\beta}, \Lie G_{-\alpha-2\beta}]
  \subset [\Lie U , V]
  \subset W
  \subset \Lie H + \Lie N
  \subset \Lie M + \Lie N .$$
  This contradicts the fact that $\Lie M + \Lie N$ does not contain nontrivial
hyperbolic elements.

\begin{subcase}
  Assume $P = P_{\alpha}$.
  \end{subcase}
  We may assume there exists $x \in \Lie H$, such that $\proj_{-\alpha} x \neq
0$ (otherwise, $H \subset N_G(N)$, so Subcase~\ref{so2npf-Pminl} applies).
Note that, because $U \subset \unip P$, we have $\proj_\alpha \Lie U = 0$.

\begin{subsubcase}
  Assume $\proj_{\alpha+\beta} \Lie U \neq 0$.
  \end{subsubcase}
  Choose $u \in \Lie U$, such that $\proj_{\alpha+\beta} u \neq 0$. Then
  $[x,u] \in [\Lie H, \Lie U] \subset \Lie U$,
and $\bigl[ [x,u], u \bigr]$ is a
nonzero element of~$\Lie G_{\alpha+2\beta}$, so we see that $\Lie
G_{\alpha+2\beta} \subset [\Lie U,\Lie U]$.
Because every unipotent subgroup of
$\SO(1,k)$ is abelian, we conclude that $\adg \Lie G_{\alpha+2\beta}$ acts
trivially on $\Lie G / \Lie H$, which means
  $ \Lie H \supset [\Lie G, \Lie G_{\alpha+2\beta}] $.
  This contradicts (\ref{Ufaithful}).

\begin{subsubcase}
  Assume $\proj_{\alpha+\beta} \Lie U = 0$.
  \end{subsubcase}
  We may assume, furthermore, that $\proj_{\alpha} \Lie H \neq 0$ (otherwise, by
replacing $H$ with its conjugate under the Weyl reflection corresponding to the
root~$\alpha$, we could revert to Subcase~\ref{so2npf-Pminl}). Then, because
$[\Lie H, \Lie U] \subset \Lie U$, we must have $\proj_{\beta} \Lie U = 0$.
Thus, $\Lie U = \Lie G_{\alpha+2\beta}$. From
Corollary~\fullref{so1k}{N-ad2inW}, we have
  $$ W
  = [\Lie G, \Lie G_{\alpha+2\beta}, \Lie G_{\alpha+2\beta}] + \Lie H
  = \Lie G_{\alpha+2\beta} + \Lie H
  \subset \Lie U + \Lie {\zar H}
  = \Lie {\zar H}
  ,$$
  so
  \begin{eqnarray*}
  W  \cap (\Lie G_{\beta} + \Lie G_{\alpha+\beta})
  &\subset& \Lie {\zar H} \cap (\Lie G_{\beta} + \Lie G_{\alpha+\beta})
  = (\Lie {\zar H} \cap \Lie N) \cap (\Lie G_{\beta} + \Lie G_{\alpha+\beta}) \\
  &=& \Lie U \cap (\Lie G_{\beta} + \Lie G_{\alpha+\beta})
  = \Lie G_{\alpha+2\beta} \cap (\Lie G_{\beta} + \Lie G_{\alpha+\beta})
  = 0 .
  \end{eqnarray*}
  On the other hand,  from Corollary~\fullref{so1k}{N-adVinW}, we know that $W$
contains a codimension-one subspace of
  $[\Lie G, \Lie G_{\alpha+2\beta}]$,
  so $W$ contains a codimension-one subspace of
  $\Lie G_{\beta} + \Lie G_{\alpha+\beta}$.
  This is a contradiction.

\begin{subcase}
  Assume $P = P_{\beta}$.
  \end{subcase}
  Note that, because $U \subset \unip P$, we have $\proj_\beta \Lie U = 0$.

  From Corollary~\fullref{so1k}{N-ad2inW}, we have
  \begin{eqnarray*}
  W
  &=& \Lie H + (\adg u)^2 \Lie G
  \subset \Lie H + (\Lie G_{\alpha} + \Lie G_{\alpha+\beta} + \Lie
G_{\alpha+2\beta}) \\
  &=& \Lie H + \unip \Lie P_\beta
  \subset (\Lie M + \Lie U) + \unip \Lie P_\beta
  = \Lie M + \unip \Lie P_\beta .
  \end{eqnarray*}

\begin{subsubcase} \label{alphazero}
  Assume there is some nonzero $u \in \Lie U$, such that $\proj_{\alpha} u = 0$.
  \end{subsubcase}
  Replacing $H$ by a conjugate (under $G_{-\beta}$), we may assume
$\proj_{\alpha+\beta} u \neq 0$.

Let $V' = V \cap (\Lie G_{-\alpha} + \Lie G_{-\alpha-\beta})$. Because $V'$
contains a codimension-one subspace of $\Lie G_{-\alpha} + \Lie
G_{-\alpha-\beta}$ (see Corollary~\fullref{so1k}{N-dimV}), one of the 
following two
subsubsubcases must apply.

\begin{subsubsubcase}
  Assume there exists $v \in V'$, such that $\proj_{-\alpha-\beta} v = 0$.
  \end{subsubsubcase}
  From Corollary~\fullref{so1k}{N-adVinW}, we have $[u,v] \in W$. Then, because
$[u,v]$ is a nonzero element of $\Lie G_{\beta}$, we conclude that
  $$0 \neq W \cap \Lie G_{\beta} \subset
(\Lie M + \unip \Lie P_\beta) \cap \Lie
G_{\beta} = 0 .$$
  This contradicts the fact that $M$, being compact, has no nontrivial unipotent
elements.

\begin{subsubsubcase}
  Assume $\proj_{-\alpha-\beta} V' = \Lie G_{-\alpha-\beta}$.
  \end{subsubsubcase}
   For $v \in V'$, we have $\proj_0[u,v] = [\proj_{\alpha+\beta} u,
\proj_{-\alpha-\beta} v]$. Thus, there is some $v \in V'$, such that $\proj_0
[u,v]$ is hyperbolic (and nonzero). On the other hand, from
Corollary~\fullref{so1k}{N-adVinW}, we have $[u,v] \in W = \Lie M + \unip \Lie
P_\beta$. This contradicts the fact that
$\Lie M \subset \Lie{\zar H}$ does not
contain nonzero hyperbolic elements.

\begin{subsubcase}
  Assume $\proj_{\alpha} u \neq 0$, for every nonzero $u \in \Lie U$.
  \end{subsubcase}
  Fix some nonzero $u \in \Lie U$. Because $\dim \Lie U_\alpha = 1$, 
we must have
$\dim \Lie U = 1$ (so $\Lie U = \real u$). Replacing $H$ by a conjugate
(under~$G_\beta$), we may assume $\proj_{\alpha+\beta} u = 0$. Also, we may
assume $\proj_{\alpha+2\beta} u \neq 0$ (otherwise, we could revert to
Subsubcase~\ref{alphazero} by replacing $H$ with its conjugate under the Weyl
reflection corresponding to the root~$\beta$).

  Let $\Lie T = [u, \Lie G_{-\alpha} + \Lie G_{-\alpha-2\beta}]$.  Because
$\langle \Lie G_\alpha, \Lie G_{-\alpha} \rangle$ and
  $\langle \Lie G_{\alpha+2\beta}, \Lie G_{-\alpha-2\beta} \rangle$
  centralize each other, we see that
  $\Lie T = [ \Lie G_\alpha, \Lie G_{-\alpha} ] + [ \Lie G_{\alpha+2\beta}, \Lie
G_{-\alpha-2\beta}]$
  is a two-dimensional subspace of~$\Lie G$ consisting entirely of hyperbolic
elements.
  Because $V$ contains a codimension-one subspace of $\Lie G_{-\alpha} + \Lie
G_{-\alpha-2\beta}$ (see Corollary~\fullref{so1k}{N-dimV}), and $[u,V]
\subset W$
(see Corollary~\fullref{so1k}{N-adVinW}), we see that $W$ contains a
codimension-one subspace of~$\Lie T$, so $W$ contains nontrivial hyperbolic
elements. This contradicts the fact that $W \subset \Lie M + \unip \Lie
P_\beta$ does not contain nontrivial hyperbolic elements.
  \end{proof}

\end{document}